\documentclass[11pt,reqno]{amsart}
\usepackage{amsmath, amsfonts, amsthm, mathtools, extarrows,amssymb,mathrsfs,dsfont}
\usepackage[x11names]{xcolor}
\usepackage[colorlinks=true,linkcolor=blue,citecolor=blue,urlcolor=black]{hyperref}
\usepackage{cleveref}
\usepackage{tikz-cd}
\usepackage[shortlabels]{enumitem}
\setlist[enumerate]{leftmargin=*}
\usepackage[margin=1.05in]{geometry}
\usepackage[indent,skip=3pt plus 1pt minus 1pt]{parskip}
\usepackage[style=alphabetic,backref,url=false,maxnames=10,giveninits,isbn=false]{biblatex}
\newtheorem{theorem}{Theorem}[section]
\newtheorem{lemma}[theorem]{Lemma}

\newtheorem{proposition}[theorem]{Proposition}

\theoremstyle{definition}
\newtheorem{definition}[theorem]{Definition}
\newtheorem{remark}[theorem]{Remark}

\numberwithin{equation}{section}

\usepackage[T1]{fontenc}

\setitemize[1]{itemsep=5pt}
\setenumerate[1]{itemsep=5pt}

\DeclareMathOperator{\SL}{SL}
\DeclareMathOperator{\Sp}{Sp}
\DeclareMathOperator{\SO}{SO}
\DeclareMathOperator{\GL}{GL}
\DeclareMathOperator{\Hom}{Hom}
\DeclareMathOperator{\km}{km}
\DeclareMathOperator{\CH}{CH}
\newcommand{\Q}{\mathbb Q}
\newcommand{\R}{\mathbb R}
\newcommand{\Z}{\mathbb Z}
\newcommand{\C}{\mathbb C}
\newcommand{\cA}{\mathcal A}
\newcommand{\cN}{\mathcal N}
\newcommand{\dd}{\mathrm{d}}

\title{Cohomological modularity of Noether--Lefschetz cycles on $\mathcal{A}_g(N)$}
\date{}
\subjclass[2020]{Primary 14K10, 11F23; Secondary 14C25, 11F27}
\keywords{Noether--Lefschetz cycles, Siegel modular varieties, modular forms, Kudla--Millson theta series, Borel--Serre compactification}
\author{Hanyu Yao}
\address{Shanghai Center for Mathematical Sciences, Fudan University, Shanghai 200438, China}
\email{hyyao23@m.fudan.edu.cn}
\begin{document}
\begin{abstract}
We prove that for every $g \geq 2$, the generating series of $d$-elliptic Noether--Lefschetz cycles on $\mathcal{A}_g(N)$ is a modular form of weight $2g$ valued in cohomology. 
For any class $[C] \in H_{2g-2}(\mathcal{A}_g(N),\mathbb{Q})$, we use the Borel--Serre boundary to construct a capped homology class $[\widehat C]$ in an orthogonal arithmetic quotient  whose pairing with the Kudla--Millson theta series equals the pairing of the generating series of $d$-elliptic Noether--Lefschetz cycles with $[C]$.
\end{abstract}

\maketitle

\section{Introduction}\label{sec:introduction}
For integers $d,g\geq1$, let $\widetilde{\mathsf{NL}}_{g,d}$ be the moduli space of morphisms of abelian varieties $h:E\to A$, where $E$ is an elliptic curve, $A$ is a principally polarized abelian variety of dimension $g$ with polarization $\Theta_A$, and $\deg(h^*\Theta_A)=d$.
The $d$-elliptic Noether--Lefschetz cycle is the class associated to the forgetful morphism $\epsilon:\widetilde{\mathsf{NL}}_{g,d}\to\cA_g$.
In \cite[Conjecture 1]{GreerLian2025Level}, Greer and Lian conjecture that the generating series of the $d$-elliptic Noether--Lefschetz cycle classes is a modular form of weight $2g$ with values in $\CH^{g-1}(\cA_g)_\Q$ if $g\geq2$.
We prove this conjecture in cohomology level.
Write $\alpha_d\in H^{2g-2}(\cA_g,\Q)$ for the cohomology class of the $d$-elliptic Noether--Lefschetz cycle and $\lambda_j=c_j(\mathbb E_g)$ for the Chern classes of the Hodge bundle. 
Put $q=e^{2\pi i\tau}$, and let $M_{2g}(\SL_2(\Z),\Q)$ denote the space of weight-$2g$ modular forms of level one with rational Fourier coefficients.

\begin{theorem}\label{thm:levelone}
For $g\geq2$, we have $$ \frac{(-1)^g}{24}\lambda_{g-1}+\sum_{d\geq1}\alpha_d q^d \in H^{2g-2}(\cA_g,\Q)\otimes M_{2g}(\SL_2(\Z),\Q).$$
\end{theorem}

The Noether--Lefschetz locus in $\cA_g$ is the locus where the Picard number exceeds one.
Debarre and Laszlo \cite{DebarreLaszlo1990} show that each irreducible component of the Noether--Lefschetz locus is either the locus of abelian varieties containing an abelian subvariety of fixed dimension $1\leq k\leq g/2$ on which the induced polarization has fixed type, or an irreducible component of a Shimura--Hilbert--Blumenthal locus arising from real multiplication.
As $d$ varies, the $d$-elliptic loci studied here exhaust the locus of abelian varieties containing an elliptic subvariety.
Indeed, the image $B=h(E)$ of a nonzero homomorphism $h:E\to A$ is an elliptic subvariety, and $\deg(h^*\Theta_A)=\deg(E\to B)\deg(\Theta_A|_B)$.
Conversely, the inclusion $B\hookrightarrow A$ realizes the degree $e=\deg(\Theta_A|_B)$, and precomposing it with isogenies onto $B$ realizes every positive multiple of $e$.
The loci obtained by fixing $e$ have codimension $g-1$ in $\cA_g$ and, for $g\geq3$, are precisely the irreducible components of maximal dimension in the Noether--Lefschetz locus.
In what follows, we refer to the $d$-elliptic Noether--Lefschetz locus simply as the Noether--Lefschetz locus.

Theorem \ref{thm:levelone} also implies the cohomological form of Greer and Lian's conjecture on $\overline{\mathcal M}_{1,1}\times\cA_g$ formulated in \cite[Conjecture 2]{GreerLian2025Level}, as we show in Section \ref{sec:modularity}.
Pulling back the level-one series along the pointed Torelli map $\mathcal M_{g,1}^{\mathrm{ct}}\to\cA_g$ yields a modular form with coefficients in $H^{2g-2}(\mathcal M_{g,1}^{\mathrm{ct}},\Q)$.
By Greer and Lian's comparison in \cite[Theorems 1.1 and 1.2]{GreerLian2025Torelli}, its coefficient of $q^d$ for $d>0$ is the cohomology class obtained by pushing forward the virtual fundamental class of degree-$d$ stable maps $(C,p)\to(E,0)$ with $C$ of compact type and $E$ varying among elliptic curves.

\subsection*{Refinement with level structure}

We refine Theorem \ref{thm:levelone} on $\cA_g(N)$, the moduli stack with full symplectic level-$N$ structure relative to $\zeta_N=e^{2\pi i/N}$.
For $d>0$ and $b\in A_N=\Hom((\Z/N\Z)^2,(\Z/N\Z)^{2g})$, let $\alpha_{d,b}(N)\in H^{2g-2}(\cA_g(N),\Q)$ be the Noether--Lefschetz cycle class obtained by requiring the map induced by $h$ on $N$-torsion to equal $b$ in the level structures on $E$ and $A$.
Set $\alpha_{0,b}(N)=\delta_{b,0}\,m_N\pi_N^*\alpha_0$, where $\alpha_0=(-1)^g\lambda_{g-1}/24$, $m_N=[\SL_2(\Z):\Gamma(N)]$, and $\pi_N:\cA_g(N)\to\cA_g$ forgets the level structure.

Let $\rho_{A_N}$ be the Weil representation of the finite quadratic module $A_N$.

\begin{theorem}\label{thm:modularity}
For $g\geq2$ and $N\geq1$, we have $$ F_N(\tau)=\sum_{b\in A_N}\sum_{d\geq0}\alpha_{d,b}(N)q^{d/N}\mathbf e_b \in H^{2g-2}(\cA_g(N),\Q)\otimes_{\Q}M_{2g}(\rho_{A_N},\Q). $$
In particular, for each $b\in A_N$, the component $F_{N,b}(\tau)=\sum_{d\geq0}\alpha_{d,b}(N)q^{d/N}$ is a modular form of weight $2g$ for $\Gamma(N)$.
\end{theorem}

\subsection*{Relation to previous work}

Iribar L\'opez showed in \cite[Corollary 4]{IribarLopez2026} that the tautological projection constructed by Canning, Molcho, Oprea, and Pandharipande in \cite{CanningMolchoOpreaPandharipande2025} sends the level-one Chow-valued series to $\frac{(-1)^g}{24}\lambda_{g-1}E_{2g}(\tau)$, where $E_{2g}$ is the normalized Eisenstein series.
For $g\leq5$, all algebraic cycle classes on $\cA_g$ are tautological in cohomology by \cite[Proposition 7]{CanningOpreaPandharipande2025}, so this formula implies Theorem \ref{thm:levelone} in that range, as Greer and Tayou observed in \cite[\S3.1]{GreerTayou2026}.
For $g=6$, the coefficient $\alpha_1=[\cA_1\times\cA_5]$ is already non-tautological in $H^{10}(\cA_6,\Q)$ according to \cite[\S1.6]{CanningOpreaPandharipande2025}.

In \cite[Theorem 1.2]{GreerLian2025Level}, Greer and Lian proved modularity of the generating series of $d$-elliptic Noether--Lefschetz classes on $\cA_1(N)\times\cA_g(N)$ for $b\in A_N$ preserving the symplectic pairings.
At level one these product classes vanish in cohomology by \cite[Proposition 3.4]{GreerLian2025Level}.
The projection to $\cA_g(N)$ is nonproper, so one cannot obtain the series with coefficients $\alpha_{d,b}(N)$ by pushing forward the product classes in ordinary cohomology.

\begin{remark}
Other proofs of cohomological modularity have been announced independently by Garc\'ia, Greer, and Iribar L\'opez and by Cheek, Tamura, and Xia.
\end{remark}

\subsection*{Strategy of the proof}

On orthogonal arithmetic manifolds, positive-norm lattice vectors define totally geodesic cycles called \emph{special cycles}, which are higher-dimensional analogues of the divisors on Hilbert modular surfaces studied by Hirzebruch and Zagier in \cite{HirzebruchZagier1976}.
Kudla and Millson proved in \cite{KudlaMillson1990} that the generating series of these special cycles is modular in cohomology.
In \cite[Proposition 5.3]{GreerLian2025Level}, Greer and Lian realized the $d$-elliptic Noether--Lefschetz cycle classes on $\cA_1(N)\times\cA_g(N)$ as pullbacks of special cycle classes under a map to an orthogonal arithmetic quotient of signature $(2g,2g)$.

Let $C$ be a rational singular $(2g-2)$-cycle on $\cA_g(N)$.
We remove small cusp neighborhoods from $\cA_1(N)$ and map the product of the remaining compact real surface with $C$ into the orthogonal arithmetic quotient.
At each cusp, the boundary is obtained by pushing a fixed cycle in a Borel--Serre face into the interior.
Assume that $N\geq3$,  we show that its class vanishes in the homology of the face.
Inspired by Funke and Millson \cite[\S4.1]{FunkeMillson2014}, we cap off the boundary to obtain a homology class $[\widehat C]$ independent of the size of the cusp neighborhoods.

We pair $[\widehat C]$ with the Kudla--Millson series and show that the pairing equals $\langle F_N(\tau),[C]\rangle$ for all $C$, giving modularity when $N\geq3$. When $N=1$ or $2$, one instead pulls the level-$N$ Kudla--Millson series back to a torsion-free level, pairs it against capped cycles, and descends along the finite cover of $\cA_g(N)$.

\noindent\textbf{Acknowledgements.}
The author is grateful to his advisor Zhiyuan Li for his constant support and many valuable comments on this work.
The author also thanks Timothy Cheek, Ryan Tamura and Jiacheng Xia as well as Luis Garc\'{i}a, Fran\c{c}ois Greer and Aitor Iribar L\'{o}pez, for their collegiality in coordinating the release of our preprints.
The author is supported by NSFC Grant No.~12121001.
\subsection*{AI disclosure}
Beginning in early 2026, the author intended to approach the problem by studying the Borel--Serre compactification, and worked out how to map the product of the Borel--Serre boundary of $\cA_1$ and $\cA_g$ into the corresponding Borel--Serre boundary of the orthogonal arithmetic manifold.
However, he did not see how to use this boundary map at that stage.

Then the author used an LLM to assist in developing the proof. The main difficulty of this project: the capped-cycle construction in Subsection \ref{subsec:capping} was provided by the LLM after the author asked it to directly construct a homology class $[\widehat C]$ on this arithmetic manifold whose pairing with the Kudla--Millson theta series equals the pairing of the generating series of Noether--Lefschetz cycles with $[C]\in H_{2g-2}(\cA_g(N),\Q)$.
The LLM also assisted with the constant-term calculation in Subsection \ref{subsec:constant}, and was used as a writing aid to improve the language and readability of the manuscript.
The author takes full responsibility for the mathematical content of the paper.

\section{Noether--Lefschetz loci and orthogonal special cycles}\label{sec:geometry}

\subsection{Noether--Lefschetz cycles with level structure}\label{subsec:moduli}
Let $\cA_g(N)$ denote the moduli stack of $g$-dimensional principally polarized abelian varieties $(A,\Theta_A)$ with full symplectic level-$N$ structure $\iota_A:A[N]\xrightarrow{\sim}(\Z/N\Z)^{2g}$. For $g=1$, we simply write $(E,\iota_E)$. Set $X=\cA_1(N)\times\cA_g(N)$, with projections $p_1$ and $p_2$.

\begin{definition}
Given $d>0$ and $b\in A_N=\Hom((\Z/N\Z)^2,(\Z/N\Z)^{2g})$, let $\widetilde{\mathsf{NL}}_{g,d}^{b}(N)$ denote the moduli stack over $X$ parametrizing homomorphisms $h:E\to A$ satisfying $\deg(h^*\Theta_A)=d$ and $\iota_A\circ h[N]=b\circ\iota_E$. Isomorphisms $h\to h'$ are pairs $(f_E,f_A)$ of polarized, level-preserving isomorphisms with $f_Ah=h'f_E$.
\end{definition}

Whereas Greer and Lian assume $b$ is symplectic in \cite[Definition 3.2]{GreerLian2025Level}, we carry out their construction for arbitrary $b\in A_N$. Even so, $\widetilde{\mathsf{NL}}_{g,d}^{b}(N)$ is empty unless $b$ pulls the standard alternating form on $(\Z/N\Z)^{2g}$ back to $d$ times that on $(\Z/N\Z)^2$, which follows on reducing the degree identity in \cite[Lemma 4.2]{GreerLian2025Level} modulo $N$.

Adding level structures on $E$ and $A$ defines a finite \'etale cover of the level-one stack $\widetilde{\mathsf{NL}}_{g,d}:=\widetilde{\mathsf{NL}}_{g,d}^{\,0}(1)$.
The matrix of $h[N]$ is locally constant on this cover, so $\widetilde{\mathsf{NL}}_{g,d}^{b}(N)$ is a union of its connected components.
It follows from the level-one case of \cite[\S3]{GreerLian2025Level} that each nonempty $\widetilde{\mathsf{NL}}_{g,d}^{b}(N)$ is smooth of dimension $1+\binom{g}{2}$.
Let $\mu^b(N)$ and $\epsilon^b(N)$ denote the forgetful maps to $\cA_1(N)$ and $\cA_g(N)$, respectively.
The stack $\widetilde{\mathsf{NL}}_{g,d}^{b}(N)$ is finite over the level-one stack, which is proper over $\cA_g$ by \cite[Lemma 3.4]{GreerLian2025Torelli}.
As $\cA_g(N)$ and $X$ are separated over $\cA_g$, both $\epsilon^b(N)$ and $(\mu^b(N),\epsilon^b(N))$ are proper.

\begin{definition}
For $d>0$ and $b\in A_N$, the Noether--Lefschetz locus for $b$ is the image of $\epsilon^b(N)$ in $\cA_g(N)$.
Let $\alpha_{d,b}(N)\in H^{2g-2}(\cA_g(N),\Q)$ denote the cycle class associated to this morphism.
\end{definition}
For $d=0$, set $\alpha_{0,b}(N)=\delta_{b,0}\,m_N\pi_N^*\alpha_0$, where $\alpha_0=(-1)^g\lambda_{g-1}/24$.

\begin{proposition}\label{prop:levels}
Assume that $N\mid M$ and denote by $\pi_{M,N}:\cA_g(M)\to\cA_g(N)$ the forgetful morphism.
Then for every $d\geq0$ and $b\in A_N$, we have $$ \sum_{\substack{b'\in A_M\\ b'\bmod N=b}}\alpha_{d,b'}(M)=\frac{m_M}{m_N}\,\pi_{M,N}^*\alpha_{d,b}(N). $$
\end{proposition}
\begin{proof}
For $d=0$, both sides are equal to $\delta_{b,0}\,m_M\pi_M^*\alpha_0$ by definition.
When $d>0$, let $$p:Y=\widetilde{\mathsf{NL}}_{g,d}^{b}(N)\times_{\cA_g(N)}\cA_g(M)\to\cA_g(M)$$ be the projection onto the second factor.
By \'etale base change along $\pi_{M,N}$, the cycle class associated to $p$ is $\pi_{M,N}^*\alpha_{d,b}(N)$.
Pull back $\cA_1(M)\to\cA_1(N)$ along $Y\to\widetilde{\mathsf{NL}}_{g,d}^{b}(N)\xrightarrow{\mu^b(N)}\cA_1(N)$ to get a map $u:\widetilde Y\to Y$.
This cover $u$ is finite \'etale of degree $m_M/m_N$, so $p_*u_*(1)=(m_M/m_N)\,p_*(1)=(m_M/m_N)\,\pi_{M,N}^*\alpha_{d,b}(N)$.
On $\widetilde Y$, the level-$M$ structures on $E$ and $A$ give a matrix for $h[M]$ that is locally constant and congruent to $b$ modulo $N$.
The substack of $\widetilde Y$ where this matrix equals a given lift $b'\in A_M$ of $b$ is $\widetilde{\mathsf{NL}}_{g,d}^{\,b'}(M)$, and on it $p\circ u$ is just $\epsilon^{b'}(M)$.
Because $\widetilde Y$ is the disjoint union of these substacks, one has $p_*u_*(1)=\sum_{b'\bmod N=b}\alpha_{d,b'}(M)$.
\end{proof}

\subsection{The tensor map and the pullback of special cycles}\label{subsec:tensor}

Consider the tensor product $L=\Lambda\otimes\Lambda'$ of the standard symplectic lattices $(\Lambda,\omega_1)$ and $(\Lambda',\omega_g)$, where we fix the symplectic bases $(e,f)$ of $\Lambda$ and $(e'_1,\ldots,e'_g,f'_1,\ldots,f'_g)$ of $\Lambda'$.
Under the product pairing $(\cdot,\cdot)=\omega_1\otimes\omega_g$ and quadratic form $Q(v)=(v,v)/2$, $L$ is an even unimodular lattice of signature $(2g,2g)$.
Write each $v\in L$ uniquely as $v=e\otimes u+f\otimes w$ with $u,w\in\Lambda'$.
Then $Q(v)=\omega_g(u,w)$, and we define $F_v:\Lambda\to\Lambda'$ by $F_v(x)=\omega_1(x,e)u+\omega_1(x,f)w$.
Through $\omega_1$, the map $v\mapsto F_v$ is an isomorphism $L\cong\Hom_\Z(\Lambda,\Lambda')$.
Let $L(N)=(L,NQ)$. 
Then we have $L/NL\cong A_N$.

At $(z,Z)\in\mathfrak H\times\mathfrak H_g$, the maps $(x,y)\mapsto x-zy$ and $(x,y)\mapsto x-Zy$ define complex structures $J_z$ on $\Lambda_\R$ and $J_Z$ on $\Lambda'_\R$, with Riemann forms $-\omega_1$ and $-\omega_g$.
For the level structures, we endow $\Lambda/N\Lambda$ and $\Lambda'/N\Lambda'$ with the induced pairings.
Write $D$ for the Grassmannian of negative $2g$-planes in $V$, with tautological bundle $\cN$.
Define the tensor map $\widetilde\phi:\mathfrak H\times\mathfrak H_g\to D$ by $\widetilde\phi(z,Z)=\ker(1+J_z\otimes J_Z)$.
The base plane $N_0=\widetilde\phi(i,iI_g)$ has an orthonormal basis $\nu_1,\ldots,\nu_{2g}$ for $-(\cdot,\cdot)$ given by $$ \nu_{2i-1}=\frac{e\otimes e'_i-f\otimes f'_i}{\sqrt2}, \qquad \nu_{2i}=-\frac{e\otimes f'_i+f\otimes e'_i}{\sqrt2} \quad(1\leq i\leq g).$$
We orient $\cN$ by pulling back the orientation of the ordered basis $(\nu_1,\ldots,\nu_{2g})$ along the orthogonal projections $P\to N_0$ for $P\in D$.
On $\widetilde\phi^*\cN$, this orientation coincides with the complex orientation from $1\otimes J_Z$ due to the relation $\nu_{2i}=(1\otimes J_{iI_g})\nu_{2i-1}$ and the connectedness of $\mathfrak H\times\mathfrak H_g$.
Because $J_z$ and $J_Z$ change by conjugation under the symplectic actions, $\widetilde\phi$ is equivariant for the tensor product action.
For a smooth choice of $s_Z\in\Sp_{2g}(\R)$ with $s_Z\cdot iI_g=Z$, we have
\begin{equation}\label{eq:periodmatrix}
 \widetilde\phi(a\cdot i,Z)=(a\otimes s_Z)N_0
 \qquad(a\in\SL_2(\R)).
\end{equation}

Denote the principal congruence subgroups of level $N$ by $\Gamma(N)\subset\SL_2(\Z)$ and $\Gamma_g(N)\subset\Sp_{2g}(\Z)$.
We also introduce $\Gamma_L(N)=\{\gamma\in\SO(L)\cap\SO(V)_0\mid\gamma\equiv1\pmod N\}$ and $\km(L(N))=[\Gamma_L(N)\backslash D]$.
The tensor product homomorphism sends $\Gamma(N)\times\Gamma_g(N)$ into $\Gamma_L(N)$, so $\widetilde\phi$ descends to $\phi:X^{\mathrm{an}}\to\km(L(N))$.
The group $\SO(V)_0$ preserves the orientation of $\cN$, and we use the same notation for the descended bundle on $\km(L(N))$.
For $N\geq3$, the groups $\Gamma(N)$, $\Gamma_g(N)$, and $\Gamma_L(N)$ are torsion free.

Set $V=L_\R$.
For $v\in V$ with $Q(v)>0$, put $D_v=\{P\in D\mid P\subset v^\perp\}$.
At $P\in D_v$, the orientation of $P$ and the basis $v$ of $\R v$ determine the orientation of the normal space $\Hom(P,\R v)$.

\begin{definition}\label{def:specialcycles}
Put $\Gamma=\Gamma_L(N)$.
For $v\in L$ with $Q(v)>0$, let $\Gamma_v$ be the stabilizer of $v$ in $\Gamma$.
The natural proper map $i_v:[\Gamma_v\backslash D_v]\to\km(L(N))$ defines a locally finite cycle $C_v$.
For $d>0$ and $b\in A_N$, define the special cycle $$ C_d^b(N)=\sum_{\substack{[v]\in\Gamma\backslash L\\ Q(v)=d,\;v\bmod NL=b}} C_v. $$
We denote the Poincar\'e dual of its fundamental class by $[C_d^b(N)]\in H^{2g}(\km(L(N)),\Q)$.
\end{definition}

For $v\in V$ with $Q(v)>0$, set $X_v=\widetilde\phi^{-1}(D_v)$.

At fixed marked periods $(z,Z)$, the homomorphisms $h:E\to A$ with $\deg(h^*\Theta_A)=d>0$ correspond bijectively to the vectors $v\in L$ with $Q(v)=d$ and $F_vJ_z=J_ZF_v$.
We have $F_{(J_z\otimes J_Z)v}=J_ZF_vJ_z^{-1}=-J_ZF_vJ_z$ because $J_z$ preserves $\omega_1$.
Note that the orthogonal complement of $\widetilde\phi(z,Z)$ is the $+1$-eigenspace of $J_z\otimes J_Z$, so $F_vJ_z=J_ZF_v$ is equivalent to $(z,Z)\in X_v$.
In the chosen level structures, $h[N]$ has matrix $F_v\bmod N$.

\begin{proposition}\label{prop:transverse}
Suppose that $v\in V$ has $Q(v)>0$.
Then $\widetilde\phi$ intersects $D_v$ transversely.
Its inverse image $X_v$ is a complex submanifold of codimension $g$, and the induced orientation on  normal bundle agrees with the complex orientation.
\end{proposition}

\begin{proof}
Equivariance under $\Sp_{2g}(\R)$ lets us assume $v=\sqrt{Q(v)}(e\otimes e'_1+f\otimes f'_1)$.
Abbreviate $P=\widetilde\phi(z,Z)$ and $J=(1\otimes J_Z)|_P$.
The $(-i)$-eigenspace of $J$ in $P_\C$ has a basis consisting of $t_j=(ze+f)\otimes\bigl(\sum_{k=1}^g Z_{kj}e'_k+f'_j\bigr)$ for $1\leq j\leq g$.
As $v$ is real, a point $(z,Z)\in\mathfrak H\times\mathfrak H_g$ lies in $X_v$ precisely when the $g$ holomorphic functions $\xi_j(z,Z)=(v,t_j(z,Z))=\sqrt{Q(v)}(\delta_{1j}+zZ_{1j})$ all vanish there.
Write $\xi=(\xi_1,\ldots,\xi_g)$.
Here $\dd\xi$ is surjective, because $z\ne0$ and a symmetric variation of $Z$ can have any first row.
In particular, $X_v$ is a complex submanifold of codimension $g$.

At a point $(z,Z)\in X_v$, set $t_j=a_j+iJa_j$, so that $a_1,Ja_1,\ldots,a_g,Ja_g$ is a positively oriented basis of $P$.
For $B\in\Hom(P,\R v)$, write $Ba_j=x_jv$ and $BJa_j=y_jv$.
Then we have $(v,Bt_j)=(x_j+iy_j)(v,v)$ with $(v,v)>0$. 
So the map $B\mapsto((v,Bt_j))_{j=1}^g$ is an orientation-preserving isomorphism $\Hom(P,\R v)\xrightarrow{\sim}\C^g$.
If $\delta$ is a tangent vector at $(z,Z)$, then $(\dd\widetilde\phi(\delta))t_j$ is the $P^\perp$-component of $\dd t_j(\delta)$.
Pairing with $v\in P^\perp$ shows that $\dd\xi$ represents the normal differential of $\widetilde\phi$ in these coordinates.
The map $\dd\xi$ is both surjective and $\C$-linear, so $\widetilde\phi$ is transverse to $D_v$ and the induced normal orientation agrees with the complex one.
\end{proof}

Fix $d>0$ and $b\in A_N$, and let $v$ range over the vectors in $L$ with $Q(v)=d$ and $v\bmod NL=b$.
Over $\mathfrak H\times\mathfrak H_g$, the period coordinates trivialize the integral homology local systems with lattices $\Lambda$ and $\Lambda'$, so the vector $v$ associated with a family of homomorphisms is locally constant.
Conversely, $F_v$ defines a homomorphism of the universal complex tori over $X_v$.
Thus $\coprod_v X_v$ parametrizes the marked homomorphisms of degree $d$ inducing $b$ on $N$-torsion.
A change of marking $(\gamma_1,\gamma_g)\in\Gamma(N)\times\Gamma_g(N)$ sends $F_v$ to $F_{(\gamma_1\otimes\gamma_g)v}=\gamma_gF_v\gamma_1^{-1}$, which has the same reduction modulo $N$.

For each $\Gamma_L(N)$-orbit $[v]$, the pullback of $i_v$ along $D\to\km(L(N))$ is the disjoint union of the inclusions $D_w\hookrightarrow D$ for $w\in\Gamma_L(N)v$.
Since $\widetilde\phi^{-1}(D_w)=X_w$, taking the quotient of the marked pullback by $\Gamma(N)\times\Gamma_g(N)$ gives the following Cartesian square, with right vertical map induced by the $i_v$.
\begin{equation}\label{eq:moduliquotient}
\begin{tikzcd}[column sep=large]
 \widetilde{\mathsf{NL}}_{g,d}^{b}(N)^{\mathrm{an}}
   \arrow[r] \arrow[d,"{(\mu^b(N),\epsilon^b(N))}"']
   & \displaystyle\coprod_{[v]}[\Gamma_v\backslash D_v] \arrow[d] \\
 X^{\mathrm{an}} \arrow[r,"\phi"'] & \km(L(N)).
\end{tikzcd}
\end{equation}
By Proposition \ref{prop:transverse}, $\widetilde\phi$ is transverse to each $D_v$, and the induced orientation of the normal bundle of $X_v$ agrees with its complex orientation.
The cycle class of $(\mu^b(N),\epsilon^b(N))_*[\widetilde{\mathsf{NL}}_{g,d}^{b}(N)]$ is therefore $\phi^*[C_d^b(N)]$ in $H^{2g}(X,\Q)$.

\subsection{The Kudla--Millson theorem}\label{subsec:KM}

Under multiplication by $N$, the discriminant module of $L(N)$ corresponds to $A_N=L/NL$ with quadratic form $q_N(\bar v)=Q(v)/N\bmod\Z$.
For the Weil representation $\rho_{A_N}$, we follow the normalization of \cite[(2.2)--(2.3)]{BruinierFunke2004}.

\begin{theorem}\label{thm:KM}
Suppose $N\mid M$, and let $e(\cN)$ be the Euler class of the negative tautological bundle at level $M$.
Then the series $$ \mathcal Z_{M,N}(\tau)=e(\cN)\mathbf e_0+\sum_{d\geq1}\sum_{b'\in A_M}[C_d^{b'}(M)]\,q^{d/N}\mathbf e_{b'\bmod N}$$ lies in $H^{2g}(\km(L(M)),\C)\otimes M_{2g}(\rho_{A_N})$.
Whenever $M\mid M'$, the natural covering $p:\km(L(M'))\to\km(L(M))$ satisfies $p^*\mathcal Z_{M,N}=\mathcal Z_{M',N}$.
\end{theorem}

\begin{proof}
Under $p$, the cycle indexed by a $\Gamma_L(M)$-orbit pulls back to the sum of the cycles indexed by the $\Gamma_L(M')$-orbits it contains.
Together with naturality of the Euler class, this gives the asserted pullback identity.
Since $p_*p^*=(\deg p)\operatorname{id}$, it suffices to prove modularity for $M\geq3$.
At such a level, apply the Kudla--Millson theorem \cite[Theorem 2]{KudlaMillson1990} in the form of \cite[Theorem 5.1]{BruinierFunke2004} to $L(N)$ and the torsion-free subgroup $\Gamma_L(M)$ of its discriminant kernel.
The vector $v/N$ defines the same cycle as $v$ and has Fourier index $Q(v)/N$, so grouping the vectors $v$ by their residues modulo $ML$ gives the displayed series.
\end{proof}

Write $\mathcal Z_N=\mathcal Z_{N,N}$.
For $N\geq3$, evaluating $\mathcal Z_N$ on a class $[\Xi]\in H_{2g}(\km(L(N)),\Q)$ gives the modular form
\begin{equation}\label{eq:KMpairing}
 \langle\mathcal Z_N(\tau),[\Xi]\rangle=\langle e(\cN),[\Xi]\rangle\mathbf e_0+\sum_{b\in A_N}\sum_{d\geq1}\langle[C_d^b(N)],[\Xi]\rangle q^{d/N}\mathbf e_b.
\end{equation}

\section{Compact cycles in the orthogonal quotient}\label{sec:capping}

Assume $N\geq3$, put $n=2g\geq4$, and let $C$ be a rational singular $(n-2)$-cycle on $\cA_g(N)$.
Inspired by the capping method of Funke and Millson \cite[\S4.1]{FunkeMillson2014}, we construct from $C$ a compact $n$-cycle in $\km(L(N))$ with homology class $[\widehat C]$.
In Section \ref{sec:coefficients}, we show that pairing the Kudla--Millson series with $[\widehat C]$ gives $\langle F_N(\tau),[C]\rangle$.

\subsection{Cusp coordinates and the Borel--Serre boundary}\label{subsec:cuspface}
We now truncate $\cA_1(N)$ near its cusps and construct the associated Borel--Serre boundary face of $\km(L(N))$ together with its maps into the interior.

For each cusp $c$ of $\cA_1(N)$, choose $\gamma_c\in\SL_2(\Z)$ such that $[\gamma_c\cdot\infty]=c$, and taking $\gamma_\infty=1$.
Every cusp of $\Gamma(N)$ has width $N$, so in the coordinate $\gamma_c^{-1}\cdot z=x+iy$, the Borel--Serre boundary circle at $c$ is parametrized by $x\bmod N$.
Set $\iota_{c,y}(x):=[\gamma_c\cdot(x+iy)]$ for $x\in\R/N\Z$ and $y>0$. 
One can choose $T_0>0$ sufficiently large that the maps $(x,y)\mapsto\iota_{c,y}(x)$ embed $(\R/N\Z)\times[T_0,\infty)$ into $\cA_1(N)$ with pairwise disjoint images as $c$ varies.
For $T\geq T_0$, let $U_c(T)=\{\iota_{c,y}(x)\mid x\in\R/N\Z,\ y>T\}$ and $$ \cA_1(N)_{\leq T}=\cA_1(N)\setminus\bigcup_c U_c(T).$$
The cusp coordinate $q_c(\iota_{c,y}(x))=\exp(2\pi i(x+iy)/N)$ identifies $U_c(T)$ with the punctured disk $0<|q_c|<e^{-2\pi T/N}$.
Removing these disjoint open disks from the compactified modular curve, we get a compact smooth real surface $\cA_1(N)_{\leq T}$ whose boundary is $\coprod_c\iota_{c,T}(\R/N\Z)$ oriented by $-\dd x$.

At the standard cusp $\infty$, the translation $z\mapsto z+N$ corresponds to $T_N=\left(\begin{smallmatrix}1&N\\0&1\end{smallmatrix}\right)$.
The operator $T_N\otimes\mathrm{id}_{\Lambda'}$ acts on $L$ via $e\otimes u+f\otimes w\mapsto e\otimes(u+Nw)+f\otimes w$ for $u,w\in\Lambda'$.
Its fixed sublattice $I=e\otimes\Lambda'$ is isotropic of rank $=\frac12\operatorname{rank}L$.
Let $P_I\subset\SO(L_\Q)$ be the rational parabolic subgroup stabilizing $I_\Q$.

By unimodularity of $\omega_g$, the tensor pairing identifies $f\otimes\Lambda'$ with $I^\vee=\Hom(I,\Z)$, sending $f\otimes w$ to the functional $e\otimes u\mapsto\omega_g(u,w)$ for $u,w\in\Lambda'$.
Under the resulting identification $L=I\oplus I^\vee$, we have $$((a,\ell),(a',\ell'))=\ell'(a)+\ell(a'), \qquad Q(a,\ell)=\ell(a).$$

Identify $I$ with $\Lambda'$ by $e\otimes a\mapsto a$, and use column coordinates in the chosen symplectic basis and its dual.
For $a,w\in\Lambda'$, the tensor $e\otimes a+f\otimes w$ has coordinates $(a,Jw)$, where $J=\left(\begin{smallmatrix}0&I_g\\-I_g&0\end{smallmatrix}\right)$ is the matrix of $\omega_g$.
Declaring the basis of $I_\R$ orthonormal identifies $I_\R^\vee$ with $I_\R$.
The base plane is then $N_0=\{(a,-a)\mid a\in I_\R\}$.

We use the horospherical coordinates for $P_I$ described in \cite[\S\S2.2--2.4]{FunkeMillson2013}.
Relative to $I_\R\oplus I_\R^\vee$, the elements of $P_I(\R)$ are the matrices $\left(\begin{smallmatrix}A&B\\0&A^{-\mathsf T}\end{smallmatrix}\right)$ with $A\in\GL(I_\R)$ and $BA^{\mathsf T}$ skew-symmetric.
We view $U_I=\bigwedge^2 I_\R$ as the space of skew-symmetric maps $I_\R^\vee\to I_\R$.
For $u\in U_I$, $y>0$, and $h\in\SL(I_\R)$, set $$ n(u)=\begin{pmatrix}1&u\\0&1\end{pmatrix},\qquad a_y=\begin{pmatrix}\sqrt y\,1&0\\0&y^{-1/2}\,1\end{pmatrix}, \qquad m(h)=\begin{pmatrix}h&0\\0&h^{-\mathsf T}\end{pmatrix}. $$
Restriction to $I_\R$ has image $\GL(I_\R)$ and kernel $n(U_I)$, the real points of the unipotent radical of $P_I$.
An element of $P_I(\R)^0$ with upper blocks $A,B$ has the unique Langlands decomposition $n(BA^{\mathsf T})a_y m(h)$, where $y=(\det A)^{2/n}$ and $h=y^{-1/2}A$.

For $h\in\SL(I_\R)$, the element $m(h)$ fixes $N_0$ if and only if $h\in\SO(I_\R)$, so the boundary symmetric space attached to $P_I$ is $D_I=\SL(I_\R)/\SO(I_\R)$.
The horospherical decomposition of $D$ with respect to $P_I$ takes the form $$ U_I\times\R_{>0}\times D_I\xrightarrow{\sim}D,\qquad (u,y,h\SO(I_\R))\longmapsto n(u)a_y m(h)N_0. $$
The corresponding Borel--Serre face of $D$ is $e(P_I)=U_I\times D_I$.
Since $a_y n(v)a_y^{-1}=n(yv)$, the point with coordinates $(u,y,\xi)$ converges to $(u,\xi)\in e(P_I)$ as $y\to\infty$ for fixed $u$ and $\xi$.
Set $\Omega=J^{-1}=-J$.
The element $\left(\begin{smallmatrix}1&x\\0&1\end{smallmatrix}\right)\otimes 1$ acts as $n(x\Omega)$ in the coordinates $(a,Jw)$.
Thus \eqref{eq:periodmatrix} takes the form $$\widetilde\phi(x+iy,Z)=n(x\Omega)a_y m(s_Z)N_0 $$ in horospherical coordinates, which converges to $(x\Omega,s_Z\SO(I_\R))\in e(P_I)$ as $y\to\infty$ for fixed $x$ and $Z$.

Put $U_I(N)=N\bigwedge^2 I$ and $\Gamma_I(N)=\ker\bigl(\SL(I)\to\SL(I/NI)\bigr)$.
Each $\gamma\in\Gamma_L(N)\cap P_I(\Q)$ restricts to an automorphism $h$ of $I$ with $h\equiv \mathrm{id}_I\pmod N$, so $\det h=1$ since $N\geq3$.
The restriction map $\Gamma_L(N)\cap P_I(\Q)\to\Gamma_I(N)$ has kernel $n(U_I(N))$ and admits the section $h\mapsto m(h)$.
Thus $\Gamma_L(N)\cap P_I(\Q)=n(U_I(N))\rtimes m(\Gamma_I(N))$, where $h\in\Gamma_I(N)$ acts on $U_I(N)$ by $u\mapsto huh^{\mathsf T}$.

\begin{definition}
The open face attached to $P_I$ in the Borel--Serre compactification of $\km(L(N))$ is $$ Y_I(N):=(\Gamma_L(N)\cap P_I(\Q))\backslash e(P_I)\cong (U_I(N)\rtimes\Gamma_I(N))\backslash(U_I\times D_I).$$
For $v\in U_I(N)$ and $\gamma\in\Gamma_I(N)$, the element $(v,\gamma)$ acts by $$(v,\gamma)\cdot(u,h\SO(I_\R))= (v+\gamma u\gamma^{\mathsf T},\gamma h\SO(I_\R)).$$
\end{definition}
Projection onto $D_I$ induces a map from $Y_I(N)$ to its Levi quotient $M_I(N)=\Gamma_I(N)\backslash D_I$.
We write $[u,h]$ for the class of $(u,h\SO(I_\R))$ in $Y_I(N)$ and $[h]$ for its image in $M_I(N)$.

The quotient $\mathcal T=U_I(N)\backslash U_I$ is a compact torus because $U_I(N)$ is a full lattice in $U_I$.
In the semidirect-product description above, $Y_I(N)=\Gamma_I(N)\backslash(\mathcal T\times D_I)$, with $\gamma\in\Gamma_I(N)$ acting on $\mathcal T$ by $[u]\mapsto[\gamma u\gamma^{\mathsf T}]$.
For $N\geq3$, Minkowski's lemma implies that $\Gamma_I(N)$ is torsion-free, so $M_I(N)$ is a smooth manifold and $D_I\to M_I(N)$ is a covering map.
Thus the projection $\pi:Y_I(N)\to M_I(N)$ is a smooth locally trivial bundle with fiber $\mathcal T$.
The zero element of $\mathcal T$, fixed by $\Gamma_I(N)$, specifies a smooth section $\sigma([h])=[0,h]$.

We now use the horospherical parametrization to map the Borel--Serre face into the interior of the orthogonal quotient at finite height.
At each fixed height $y>0$, we have a smooth map $U_I\times D_I\to D$ sending $(u,h\SO(I_\R))$ to $n(u)a_y m(h)N_0$.
For $v\in U_I(N)$ and $\gamma\in\Gamma_I(N)$, the identities $m(\gamma)n(u)m(\gamma)^{-1}=n(\gamma u\gamma^{\mathsf T})$ and $m(\gamma)a_y=a_y m(\gamma)$ yield $$n(v+\gamma u\gamma^{\mathsf T})a_y m(\gamma h)N_0=n(v)m(\gamma)n(u)a_y m(h)N_0.$$
The relation above holds for every $y>0$, and $n(v)m(\gamma)$ belongs to $\Gamma_L(N)$, so the family descends to a smooth map $(0,\infty)\times Y_I(N)\to\km(L(N))$.
Its restriction at height $y$ is $$j_y:Y_I(N)\longrightarrow\km(L(N)), \qquad j_y([u,h])=[n(u)a_y m(h)N_0].$$

At the standard cusp, the $U_I$ and $D_I$ coordinates of $\widetilde\phi(x+iy,Z)$ define a smooth map $\R\times\mathfrak H_g\to U_I\times D_I$, $(x,Z)\mapsto(x\Omega,s_Z\SO(I_\R))$.
Replacing $x$ by $x+N$ leaves the image in $Y_I(N)$ unchanged because $N\Omega\in U_I(N)$.
For $\gamma\in\Gamma_g(N)\subset\Gamma_I(N)$, the identities $\gamma\Omega\gamma^{\mathsf T}=\Omega$ and $s_{\gamma Z}^{-1}\gamma s_Z\in U(g)\subset\SO(I_\R)$ ensure compatibility with the action on $\mathfrak H_g$.
The resulting smooth map
$$f:(\R/N\Z)\times\cA_g(N)\longrightarrow Y_I(N), \qquad f([x],[Z])=[x\Omega,s_Z]$$
satisfies $\phi\circ(\iota_{\infty,y}\times\mathrm{id})=j_y\circ f$ for $y>0$.
For each $[Z]$, the map $f$ sends $(\R/N\Z)\times\{[Z]\}$ into the torus fiber $\pi^{-1}([s_Z])$.

The construction at $\infty$ transports to every cusp $c$ via conjugation with $\gamma_c$.
At each $c$, the conjugate of $z\mapsto z+N$ fixes $I_c=(\gamma_c\otimes1)I$, whose associated parabolic subgroup is $(\gamma_c\otimes1)P_I(\gamma_c\otimes1)^{-1}$.
Moreover, $\gamma_c\otimes1\in\SO(V)_0$ is an integral isometry normalizing $\Gamma_L(N)$, and therefore induces a diffeomorphism $\overline\gamma_c$ of $\km(L(N))$.
Because $j_{c,y}=\overline\gamma_c\circ j_y$, we have
\begin{equation}\label{eq:allcusps}
\phi\circ(\iota_{c,y}\times\mathrm{id})=j_{c,y}\circ f.
\end{equation}
On the Borel--Serre boundary, $\gamma_c\otimes1$ identifies $Y_I(N)$ with the face attached to $I_c$.
The boundary map at $c$ is the composite of $f$ with this identification.
To construct the caps, it remains to show that $f_*([\R/N\Z]\times[C])$ vanishes in $H_{n-1}(Y_I(N),\Q)$.

\subsection{Vanishing of the boundary class}
For Hilbert modular surfaces, Funke and Millson \cite[Lemma 3.1 and Proposition 4.2(ii)]{FunkeMillson2014} show that every circle in a torus fiber bounds a rational $2$-chain in the corresponding Borel--Serre boundary component.
For the family of cusp circles defined by $f$, we use the Serre spectral sequence of $\pi$ to prove that $f_*([\R/N\Z]\times[C])=0$ for every $[C]\in H_{n-2}(\cA_g(N),\Q)$.

For the torus $\mathcal T$ defined above, put
$$W=H_1(\mathcal T,\Q)\simeq U_I(N)\otimes_\Z\Q \simeq \bigwedge\nolimits^2 I_\Q.$$
These identifications are $\Gamma_I(N)$-equivariant for the exterior-square action on $\bigwedge^2 I_\Q$.

\begin{lemma}\label{lem:vanishing}
For $n=2g\geq4$ and $N\geq3$, we have $H_{n-2}\bigl(\Gamma_I(N),\bigwedge^2 I_\Q\bigr)=0$.
\end{lemma}

\begin{proof}
Put $W_\C=W\otimes_\Q\C$ and let $F=W_\C^\vee$ be its contragredient representation.
The $\GL(I_\C)$-representation $F$ is irreducible, and its restriction to $\SL(I_\C)$ is irreducible and nontrivial.
Zariski density of $\Gamma_I(N)$ in $\SL(I_\C)$ therefore implies $F^{\Gamma_I(N)}=F^{\SL(I_\C)}=0$.

For $N\geq3$, the principal congruence subgroup of $\GL(I)$ of level $N$ is precisely $\Gamma_I(N)$, since $h\equiv1\pmod N$ and $\det h=\pm1$ imply $\det h=1$.
Li and Sun's comparison theorem \cite[Theorem 1.8 and Example 1.10]{LiSun2019} applies to $\GL_n(\R)$ in degrees $i<n-1$.
Here $F^{\Gamma_I(N)}=0$, and the theorem yields $H^{n-2}(\Gamma_I(N),F)=0$.

Duality between group homology and cohomology over $\C$ gives $$ H_{n-2}(\Gamma_I(N),W_\C)^\vee \simeq H^{n-2}(\Gamma_I(N),F)=0.$$
Group homology commutes with extension of scalars, and hence $H_{n-2}(\Gamma_I(N),W)=0$.
\end{proof}

\begin{proposition}\label{prop:boundary}
The class $f_*([\R/N\Z]\times[C])$ vanishes in $H_{n-1}(Y_I(N),\Q)$.
\end{proposition}

\begin{proof}
Put $\beta=f_*([\R/N\Z]\times[C])$ and fix an integer $r>1$.
Multiplication $u\mapsto ru$ maps $U_I(N)$ into itself and commutes with the $\Gamma_I(N)$-action, so it induces a smooth bundle map $[r]:Y_I(N)\to Y_I(N)$ over the identity on $M_I(N)$, given by $[r]([u,h])=[ru,h]$.
Let $m_r:\R/N\Z\to\R/N\Z$ be multiplication by $r$.
By the definition of $f$, we have $[r]\circ f=f\circ(m_r\times\mathrm{id})$.
Since $m_r$ has degree $r$, naturality of the homology cross product yields $[r]_*\beta=r\beta$.

The free, properly discontinuous action of $\Gamma_I(N)$ on the contractible space $D_I$ makes $M_I(N)$ a $K(\Gamma_I(N),1)$.
The fiber homology local system in degree $q$ is associated with the $\Gamma_I(N)$-module $H_q(\mathcal T,\Q)\simeq\bigwedge\nolimits^q W$.
The homology Serre spectral sequence of $\pi:Y_I(N)\to M_I(N)$ therefore has the form $$ E^2_{p,q}=H_p\bigl(\Gamma_I(N),\bigwedge\nolimits^q W\bigr) \Longrightarrow H_{p+q}(Y_I(N),\Q).$$
The map $[r]$ covers the identity on the base and acts on $H_q(\mathcal T,\Q)$ by the scalar $r^q$, which is also its action on $E^2_{p,q}$.
By naturality of the differentials, the same scalar acts on $E^s_{p,q}$ for every $s\geq2$ and hence on $E^\infty_{p,q}$.

Set $H=H_{n-1}(Y_I(N),\Q)$.
By naturality, $[r]_*$ preserves the Serre filtration $0=F_{-1}H\subseteq F_0H\subseteq\cdots\subseteq F_{n-1}H=H$, and the isomorphisms $F_pH/F_{p-1}H\simeq E^\infty_{p,n-1-p}$ are equivariant.
Consequently, $[r]_*$ acts on $F_pH/F_{p-1}H$ by multiplication by $r^{n-1-p}$.

Suppose that $\beta\ne0$ and choose the least $p$ such that $\beta\in F_pH$.
Its image $\bar\beta$ in $F_pH/F_{p-1}H$ is nonzero and satisfies $[r]_*\bar\beta=r\bar\beta$.
Comparison with the action of $[r]_*$ on the quotient yields $(r^{n-1-p}-r)\bar\beta=0$, forcing $p=n-2$ over the rationals for $r>1$.
However, $E^2_{n-2,1}=H_{n-2}(\Gamma_I(N),W)=0$ by Lemma \ref{lem:vanishing}, so $E^\infty_{n-2,1}=0$.
This contradicts $\bar\beta\ne0$ and proves $\beta=0$.
\end{proof}

\subsection{Capping the boundary}\label{subsec:capping}

The factorization \eqref{eq:allcusps} allows the boundaries of the truncated products at different heights to be represented by images of fixed cycles in $Y_I(N)$.
Adapting the capping method of Funke and Millson \cite[\S4.1]{FunkeMillson2014} and their finite-height construction \cite[\S4.2]{FunkeMillson2011Spectacle}, we choose caps for these cycles in $Y_I(N)$ and map them into $\km(L(N))$ at varying heights.

We work with smooth singular chains with rational coefficients.
Their inclusion into the ordinary singular chain complex induces an isomorphism on homology by \cite[Theorem 18.7]{Lee2013}.
Choose a smooth representative of $[C]$, again denoted by $C$.
For products of chains, we use the standard Eilenberg--Zilber shuffle product, which is associative and natural at the chain level.

Fix $T_*\geq T_0$ and choose a finite smooth triangulation of $\cA_1(N)_{\leq T_*}$ with a total ordering of its vertices.
Let $F_{T_*}$ be the sum of the ordered $2$-simplices with signs determined by the complex orientation.
The interior edges cancel in pairs in $\partial F_{T_*}$, leaving the signed edges on each boundary circle.
For each cusp $c$, let $s_c$ be the $1$-cycle on $\R/N\Z$ whose pushforward under $\iota_{c,T_*}$ is the signed sum of the boundary edges at $c$.
Then $s_c$ represents the fundamental class of $\R/N\Z$ with orientation $-\dd x$, and $\partial F_{T_*}=\sum_c(\iota_{c,T_*})_*s_c$.

For $T>T_*$, define $\kappa_{c,T}:[0,1]\times(\R/N\Z)\to\cA_1(N)$ by $\kappa_{c,T}(t,x)=\iota_{c,T_*+t(T-T_*)}(x)$.
Write $[0,1]$ for the identity singular $1$-simplex.
Extend $F_{T_*}$ over the cusp cylinders by setting 
$$ F_T=F_{T_*}+\sum_c(\kappa_{c,T})_*([0,1]\times s_c). $$
In cusp coordinates, $\kappa_{c,T}^*(\dd x\wedge\dd y)=(T-T_*)\,\dd t\wedge(-\dd x)$, so the product orientation agrees with the complex orientation.
The $c$th cylinder contributes $(\iota_{c,T})_*s_c-(\iota_{c,T_*})_*s_c$ to the boundary.
Thus, for every $T\geq T_*$, the chain $F_T$ represents the relative fundamental class of $\cA_1(N)_{\leq T}$ and satisfies $$ \partial F_T=\sum_c(\iota_{c,T})_*s_c.$$

Put $z_c=f_*(s_c\times C)$, a smooth rational $(n-1)$-cycle in $Y_I(N)$.
Using $\partial C=0$ and naturality of the cross product, we deduce from \eqref{eq:allcusps} that $$ \partial\phi_*(F_T\times C)=\phi_*((\partial F_T)\times C)=\sum_c(j_{c,T})_*z_c.$$

Recall that $\pi:Y_I(N)\to M_I(N)$ is the torus bundle projection.
\begin{lemma}\label{lem:caps}
For each cusp $c$, there exists a smooth singular $n$-chain $K_c$ in $Y_I(N)$ with rational coefficients such that $\partial K_c=z_c$ and $\int_{K_c}\pi^*\eta=0$ for every smooth $n$-form $\eta$ on $M_I(N)$.
\end{lemma}

\begin{proof}
By Proposition \ref{prop:boundary} and the comparison with smooth singular homology, choose a smooth rational $n$-chain $K'_c$ in $Y_I(N)$ with $\partial K'_c=z_c$.
We first construct a chain $R_c$ in $M_I(N)$ with the same boundary as $\pi_*K'_c$ and with zero integral against every smooth $n$-form.

Choose a closed disk $D_c$ and a diffeomorphism $i_c:\R/N\Z\to\partial D_c$.
The disk $D_c$ is contractible, so we may choose a smooth rational $2$-chain $b_c$ in $D_c$ with $\partial b_c=(i_c)_*s_c$.
Define $\widetilde f_c:D_c\times\cA_g(N)\to M_I(N)$ by $\widetilde f_c(q,[Z])=[s_Z]$, so that $\widetilde f_c\circ(i_c\times\mathrm{id})=\pi\circ f$.
Set $R_c=(\widetilde f_c)_*(b_c\times C)$.
Using $\partial C=0$, we obtain $$ \partial R_c=(\widetilde f_c)_*((i_c)_*s_c\times C) =(\pi\circ f)_*(s_c\times C)=\pi_*z_c. $$
As $\widetilde f_c$ is independent of the disk coordinate, each singular simplex in $R_c$ factors smoothly through $\Delta^{n-2}$.
Every smooth $n$-form $\eta$ on $M_I(N)$ therefore pulls back to zero on each such simplex, so $\int_{R_c}\eta=0$.

Set $K_c=K'_c-\sigma_*(\pi_*K'_c-R_c)$.
The chain $\pi_*K'_c-R_c$ is closed, so $\partial K_c=z_c$.
Moreover, $\pi\circ\sigma=\mathrm{id}$ implies $\pi_*K_c=R_c$, and hence $\int_{K_c}\pi^*\eta=\int_{R_c}\eta=0$ for every smooth $n$-form $\eta$ on $M_I(N)$.
\end{proof}

Choose the caps $K_c$ once and for all.
For $T\geq T_*$, define
\begin{equation}\label{eq:cappedcycle}
 \widehat C_T=\phi_*(F_T\times C)-\sum_c(j_{c,T})_*K_c.
\end{equation}

\begin{proposition}\label{prop:compactcycle}
For every $T\geq T_*$, the chain $\widehat C_T$ is a smooth singular $n$-cycle with rational coefficients in $\km(L(N))$.
Its homology class $[\widehat C_T]\in H_n(\km(L(N)),\Q)$ does not depend on $T$.
\end{proposition}

\begin{proof}
We have $\partial\widehat C_T=\sum_c(j_{c,T})_*(z_c-\partial K_c)=0$, so $\widehat C_T$ is a cycle.

Now fix $T>T_*$ and set $y(t)=T_*+t(T-T_*)$ for $0\leq t\leq1$.
For each cusp $c$ and each smooth singular chain $A$ in $Y_I(N)$, let $H_{c,T}(A)$ be the pushforward of $[0,1]\times A$ under the map $(t,p)\mapsto j_{c,y(t)}(p)$.
We have $\partial H_{c,T}(K_c)=(j_{c,T})_*K_c-(j_{c,T_*})_*K_c-H_{c,T}(z_c)$.
On each cusp cylinder, \eqref{eq:allcusps} shows that $\phi(\kappa_{c,T}(t,x),[Z])=j_{c,y(t)}(f(x,[Z]))$.
Since $F_T-F_{T_*}$ is the sum of these cylinders and $z_c=f_*(s_c\times C)$, naturality and associativity of the cross product imply $$ \phi_*\bigl((F_T-F_{T_*})\times C\bigr)=\sum_c H_{c,T}(z_c).$$
It follows that $\widehat C_T-\widehat C_{T_*}=-\partial\bigl(\sum_c H_{c,T}(K_c)\bigr)$.
\end{proof}

Write $[\widehat C]$ for this common homology class.
In the coefficient calculations, we may therefore vary the truncation height $T$ without affecting the pairing with the Kudla--Millson generating series.

\section{Comparison of Fourier coefficients}\label{sec:coefficients}

Assume $N\geq3$ and fix $C$ and $[\widehat C]$ as in Section \ref{sec:capping}.
We compute the positive and constant Fourier coefficients of \eqref{eq:KMpairing} for $\Xi=[\widehat C]$.

\subsection{Positive Fourier coefficients}
To compute the positive Fourier coefficients, we first establish the following analogue of \cite[Proposition 4.1]{FunkeMillson2014}.

\begin{lemma}\label{lem:avoidance}
Let $K\subset Y_I(N)$ be compact.
There is a height $T(K,d)$ such that $j_{c,y}(K)$ is disjoint from the support of $C_d^b(N)$ for every $y\geq T(K,d)$, every cusp $c$, and every $b\in A_N$.
\end{lemma}

\begin{proof}
Write $S_d=\bigcup_{b\in A_N}\operatorname{supp}C_d^b(N)$, which lifts in $D$ to $\bigcup_{v\in L,\,Q(v)=d} D_v$.
Since $\gamma_c\otimes 1$ is an isometry of $L$, it preserves this preimage and hence $\overline\gamma_c(S_d)=S_d$.
As $j_{c,y}=\overline\gamma_c\circ j_y$, it suffices to prove the assertion for $j_y$.

The quotient map $U_I\times\SL(I_\R)\to Y_I(N)$ is open and its source is locally compact, so the images of finitely many relatively compact open sets cover $K$.
Let $\widetilde K$ be the union of their closures.
By compactness, there is a constant $B_K>0$ such that $\|h^{-\mathsf T}\xi\|\leq B_K\|\xi\|$ for every $(u,h)\in\widetilde K$ and every $\xi\in I_\R^\vee$.

Suppose that $j_y(k)\in S_d$ for some $k\in K$, choose a representative $(u,h)\in\widetilde K$ of $k$, and put $G=n(u)a_y m(h)$.
Then $GN_0\in D_v$ for some $v=(a,\ell)\in L=I\oplus I^\vee$ with $Q(v)=d$.
In the coordinates of Section \ref{subsec:cuspface}, we have $N_0^\perp=\{(x,x)\mid x\in I_\R\}$, so $v\perp GN_0$ implies $v=G(x,x)$ for some $x\in I_\R$.
As $G$ preserves $Q$, we have $\|x\|^2=Q(x,x)=d$, and the second component of $G(x,x)$ is $\ell=y^{-1/2}h^{-\mathsf T}x$.
The equality $\ell(a)=d>0$ implies that $\ell$ is nonzero.
Its coordinates in the orthonormal dual basis are integers, so $$ 1\leq \|\ell\|^2=y^{-1}\|h^{-\mathsf T}x\|^2 \leq \frac{B_K^2d}{y}.$$
Thus any $T(K,d)>B_K^2d$ has the required property.
\end{proof}
\begin{proposition}\label{prop:positive}
For $d>0$ and $b\in A_N$, one has $\bigl\langle[C_d^b(N)],[\widehat C]\bigr\rangle =\bigl\langle\alpha_{d,b}(N),[C]\bigr\rangle$.
\end{proposition}
\begin{proof}
Put $S=\operatorname{supp}C_d^b(N)$.
Choose a relatively compact open neighborhood $U$ of $\operatorname{supp}C$ in $\cA_g(N)$, and write $[C]_U$ for the class of $C$ in $H_{2g-2}(U,\Q)$.
The set $K=f((\R/N\Z)\times\overline U)\cup\bigcup_c\operatorname{supp}K_c$ is compact.
By Lemma \ref{lem:avoidance}, we may choose $T\geq T_*$ such that $j_{c,y}(K)\cap S=\varnothing$ for every cusp $c$ and every $y\geq T$.
In particular, the cap images $j_{c,T}(\operatorname{supp}K_c)$ are disjoint from $S$.
Writing $M_T=\cA_1(N)_{\leq T}$, we deduce from the cusp factorization \eqref{eq:allcusps} that every point of $\cA_1(N)\times U$ whose image under $\phi$ lies in $S$ belongs to $\operatorname{int}M_T\times U$.

Let $\mathcal H_U=(\epsilon^b(N))^{-1}(U)$ and let $[\mathcal H_U]$ be its Borel--Moore fundamental class.
Since the image of $\mathcal H_U$ under $(\mu^b(N),\epsilon^b(N))$ lies in $\operatorname{int}M_T\times U$, the proper map $r=(\mu^b(N),\epsilon^b(N)):\mathcal H_U\to M_T\times U$ is the pullback of the right vertical map in \eqref{eq:moduliquotient} along $\phi|_{M_T\times U}$.
By Proposition \ref{prop:transverse}, this pullback is transverse and its induced normal orientation matches the complex orientation.
The caps $(j_{c,T})_*K_c$ lie in $\km(L(N))\setminus S$, so $\widehat C_T$ and $\phi_*(F_T\times C)$ define the same class in $H_{2g}(\km(L(N)),\km(L(N))\setminus S;\Q)$.
Since $F_T\times C$ represents $[M_T,\partial M_T]\times[C]_U$, transverse pullback implies $$ \bigl\langle[C_d^b(N)],[\widehat C]\bigr\rangle=\bigl([M_T,\partial M_T]\times[C]_U\bigr)\cdot r_*[\mathcal H_U].$$
The projection formula for the proper map $M_T\times U\to U$ then yields $$ \bigl([M_T,\partial M_T]\times[C]_U\bigr)\cdot r_*[\mathcal H_U]=[C]_U\cdot\bigl(\epsilon^b(N)|_{\mathcal H_U}\bigr)_*[\mathcal H_U]. $$
Since proper pushforward commutes with restriction to $U$, the class $\bigl(\epsilon^b(N)|_{\mathcal H_U}\bigr)_*[\mathcal H_U]$ is dual to $\alpha_{d,b}(N)|_U$.
The right-hand side therefore equals $\langle\alpha_{d,b}(N)|_U,[C]_U\rangle=\langle\alpha_{d,b}(N),[C]\rangle$.
\end{proof}

\subsection{The constant coefficient\label{subsec:constant}}

We compute the pairing of the Euler class $e(\cN)$ with the capped class $[\widehat C]$.
On $D$, a local smooth section $s$ of $\cN$ is a $V$-valued map with $s(P)\in P$.
For $X\in T_PD$, define $\nabla^{\cN}_X s$ at $P$ to be the projection of $\dd s_P(X)$ onto $P$ along $P^\perp$.
The connection $\nabla^{\cN}$ is $\SO(V)_0$-invariant and preserves the metric $-(\cdot,\cdot)$, so its Euler form descends to a form $e_{\cN}$ on $\km(L(N))$.
By Proposition \ref{prop:compactcycle} and \eqref{eq:cappedcycle},  we have
\begin{equation}\label{eq:constantpairing}
 \langle e(\cN),[\widehat C]\rangle=\int_{F_T\times C}\phi^*e_{\cN}-\sum_c\int_{K_c}j_{c,T}^*e_{\cN}
\end{equation}
for every $T\geq T_*$.

Let $\mathbb E_{r,N}$ be the Hodge bundle on $\cA_r(N)$, and let $\lambda_{j,N}$ be the $j$th Chern class of $\mathbb E_{g,N}$.
We endow $\mathbb E_{r,N}$ with the Hodge metric induced by the principal polarization, normalized so that $\dd w_1,\ldots,\dd w_r$ are orthonormal at $iI_r$.
We denote by $c_j$ the $j$th Chern form with respect to the induced Hermitian metrics.
\begin{lemma}\label{lem:tensorEuler}
For every $N\geq1$, we have
\begin{equation}\label{eq:tensorEuler}
 \phi^*e_{\cN}=c_g(\mathcal F)=p_2^*c_g(\mathbb E_{g,N}^{\vee})+p_1^*c_1(\mathbb E_{1,N}^{\vee})\wedge p_2^*c_{g-1}(\mathbb E_{g,N}^{\vee}).
\end{equation}
\end{lemma}

\begin{proof}
Over the genus-$r$ period domain, let $H_r^\pm$ be the $\pm i$-eigenbundles of the complex structure on the complexified first homology of the universal abelian variety.
At a period matrix $W$, the map $(x,y)\mapsto x-Wy$ from complexified homology to the tangent space of the abelian variety at the origin has kernel $H_r^-$.
Its restriction to $H_r^+$ identifies this smooth bundle with the dual Hodge bundle $\mathbb E_{r,N}^{\vee}$.
Write $\nabla_r$ for the Gauss--Manin connection, and for a local smooth section $s$ of $H_r^+$ let $D_rs$ be the projection of $\nabla_rs$ onto $H_r^+$ along $H_r^-$.
By the description of the Chern connection for a pure polarized variation in \cite[\S5.3]{PearlsteinPeters2019}, $D_r$ is the Chern connection of $\mathbb E_{r,N}^{\vee}$ under this identification.

As in \cite[proof of Theorem 1.2]{GreerLian2025Level}, we have $$ (\phi^*\cN)_\C=(H_1^+\otimes H_g^+)\oplus(H_1^-\otimes H_g^-)=\mathcal F\oplus\overline{\mathcal F}. $$
Here $\mathcal F$ is the $+i$-eigenbundle of $1\otimes J_Z$ on $(\phi^*\cN)_\C$.
For local sections $u$ of $H_1^+$ and $v$ of $H_g^+$, the Gauss--Manin connection on the tensor product of the complexified first homology bundles sends $u\otimes v$ to $(\nabla_1u)\otimes v+u\otimes(\nabla_gv)$.
Orthogonal projection onto $(\phi^*\cN)_\C$ annihilates $H_1^-\otimes H_g^+$ and $H_1^+\otimes H_g^-$, so $(\phi^*\nabla^{\cN})(u\otimes v)=D_1u\otimes v+u\otimes D_gv$.
Hence $\phi^*\nabla^{\cN}$ preserves $\mathcal F$ and restricts there to the tensor product Chern connection.
By complex conjugation, it also preserves $\overline{\mathcal F}$.
With the complex orientation, $\phi^*e_{\cN}$ is the top Chern form $c_g(\mathcal F)$.

Put $\ell=p_1^*c_1(\mathbb E_{1,N}^{\vee})$.
One has $$ c_g(\mathcal F)=\sum_{j=0}^g \ell^j\wedge p_2^*c_{g-j}(\mathbb E_{g,N}^{\vee}).$$
As $\ell$ is pulled back from the modular curve, we have $\ell^2=0$.
So only the $j=0,1$ terms survive.
\end{proof}

For $N\geq3$, let $e_I=(j_1\circ\sigma)^*e_{\cN}$ on $M_I(N)$, where $j_1$ is the map at height $1$ for the standard cusp.

\begin{lemma}\label{lem:Eulerboundary}
For every cusp $c$, the forms $j_{c,y}^*e_{\cN}$ converge to $\pi^*e_I$ in $C^\infty$ on compact subsets of $Y_I(N)$ as $y\to\infty$.
In particular, $\displaystyle\lim_{y\to\infty}\int_{K_c}j_{c,y}^*e_{\cN}=0$ for each cap chosen in Lemma \ref{lem:caps}.
\end{lemma}

\begin{proof}
Write $e_D$ for the lift of $e_{\cN}$ to $D$.
In the coordinates of Section \ref{subsec:cuspface}, define $\psi_t:U_I\times D_I\to D$ for $t\in\R$ by $$\psi_t(u,h\SO(I_\R))=n(tu)m(h)N_0.$$
The identity $n(u)a_y=a_y n(u/y)$ shows that the composite $U_I\times D_I\to Y_I(N)\xrightarrow{j_{c,y}}\km(L(N))$ lifts to $(\gamma_c\otimes1)a_y\psi_{1/y}$.
By invariance of $e_D$, the lift of $j_{c,y}^*e_{\cN}$ is therefore $\psi_{1/y}^*e_D$.
The map $(t,p)\mapsto\psi_t(p)$ is smooth on $\R\times(U_I\times D_I)$, so $\psi_{1/y}^*e_D$ converges to $\psi_0^*e_D$ in $C^\infty$ on compact subsets.
At $t=0$, the map $(u,h\SO(I_\R))\mapsto m(h)N_0$ lifts $j_1\circ\sigma\circ\pi$, where $\sigma([h])=[0,h]$.
The asserted convergence on $Y_I(N)$ follows by taking local lifts to $U_I\times D_I$.

The convergence is uniform on the compact support of the fixed smooth chain $K_c$, so
$$\lim_{y\to\infty}\int_{K_c}j_{c,y}^*e_{\cN}=\int_{K_c}\pi^*e_I=0.$$
The last equality follows from Lemma \ref{lem:caps} with $\eta=e_I$.
\end{proof}

\begin{proposition}\label{prop:constant}
For $N\geq3$, we have $$ \langle e(\cN),[\widehat C]\rangle=\frac{(-1)^g m_N}{24}\langle\lambda_{g-1,N},[C]\rangle. $$
\end{proposition}

\begin{proof}
The $2g$-form $c_g(\mathbb E_{g,N}^{\vee})$ pulls back to zero on each $(2g-2)$-simplex of $C$, so only the mixed term in \eqref{eq:tensorEuler} contributes to the integral over $F_T\times C$.
As $F_T$ represents the relative fundamental class of $M_T$, integration over the two factors yields
$$ \begin{aligned}
 \int_{F_T\times C}\phi^*e_{\cN}
 &=\left(\int_{M_T}c_1(\mathbb E_{1,N}^{\vee})\right) \left(\int_C c_{g-1}(\mathbb E_{g,N}^{\vee})\right)\\
 &=(-1)^g\left(\int_{M_T}c_1(\mathbb E_{1,N})\right)
   \langle\lambda_{g-1,N},[C]\rangle.
 \end{aligned} $$

For $z=x+iy$, the differential $\dd w$ spans the Hodge line of $E_z=\C/(\Z+z\Z)$ with squared Hodge norm $\|\dd w\|_z^2=\tfrac{i}{2}\int_{E_z}\dd w\wedge\dd\overline w=y$.
Hence $$ c_1(\mathbb E_{1,N})=-\frac{i}{2\pi}\partial\bar\partial\log y=\frac{1}{4\pi}\frac{\dd x\wedge\dd y}{y^2}.$$
Since $-I_2$ acts trivially on $\mathfrak H$ and lies outside $\Gamma(N)$ for $N\geq3$, the hyperbolic area of $\cA_1(N)$ is $(m_N/2)\cdot(\pi/3)$.
As the regions $M_T$ exhaust $\cA_1(N)$, one has $$ \lim_{T\to\infty}\int_{M_T}c_1(\mathbb E_{1,N})=\frac{1}{4\pi}\cdot\frac{m_N}{2}\cdot\frac{\pi}{3}=\frac{m_N}{24}.$$
In \eqref{eq:constantpairing}, the cap integrals vanish as $T\to\infty$ by Lemma \ref{lem:Eulerboundary}, so the proposition follows.
\end{proof}
\section{Modularity of the generating series}\label{sec:modularity}

\begin{proof}[Proof of Theorem \ref{thm:modularity}]
Choose a multiple $M\geq3$ of $N$, and write $\pi=\pi_{M,N}:\cA_g(M)\to\cA_g(N)$.
Fix a homology class $[C]\in H_{2g-2}(\cA_g(M),\Q)$, and let $[\widehat C]$ be the capped class constructed in Section \ref{sec:capping}.
By Proposition \ref{prop:positive}, the coefficient of $q^{d/N}\mathbf e_b$ in $\langle\mathcal Z_{M,N}(\tau),[\widehat C]\rangle$ for $d>0$ and $b\in A_N$ is $\sum_{b'\in A_M,\,b'\bmod N=b}\langle\alpha_{d,b'}(M),[C]\rangle$.
The same formula holds for $d=0$ by Proposition \ref{prop:constant}, since only $b'=0$ contributes to the constant term.
Together with Proposition \ref{prop:levels}, this yields $$ \bigl\langle\pi^*F_N(\tau),[C]\bigr\rangle=\frac{m_N}{m_M}\bigl\langle\mathcal Z_{M,N}(\tau),[\widehat C]\bigr\rangle.$$

By Theorem \ref{thm:KM}, $\langle\pi^*F_N,[C]\rangle$ belongs to $M_{2g}(\rho_{A_N})$ for every $[C]$.
Evaluating on a basis of $H_{2g-2}(\cA_g(M),\Q)$, which is finite-dimensional by \cite[\S11.1(c)]{BorelSerre1973}, we obtain $\pi^*F_N\in H^{2g-2}(\cA_g(M),\C)\otimes M_{2g}(\rho_{A_N})$.
The map $\pi$ is representable and finite \'etale, so $\pi_*\pi^*=(\deg\pi)\operatorname{id}$.
Applying $(\deg\pi)^{-1}\pi_*$ to the cohomology factor proves the vector-valued assertion for $F_N$.

Each component $F_{N,b}$ is modular of weight $2g$ for $\Gamma(N)$ because $\rho_{A_N}$ is trivial on this subgroup.

\end{proof}

\begin{remark}
Theorem \ref{thm:levelone} also implies the cohomological version of Greer and Lian's modularity conjecture on $\overline{\mathcal M}_{1,1}\times\cA_g$.
Let $\overline\alpha_d\in H^{2g}(\overline{\mathcal M}_{1,1}\times\cA_g,\Q)$ be the pushforward of $[\widetilde{\mathsf{NL}}_{g,d}]$ along $(E,A,h)\mapsto(E,A)$ for $d>0$, and put $\overline\alpha_0=c_g(\mathbb E_1^\vee\otimes\mathbb E_g^\vee)$, where $\mathbb E_1$ and $\mathbb E_g$ are pulled back from the two factors.
For the open immersion $j:\cA_1\times\cA_g\hookrightarrow\overline{\mathcal M}_{1,1}\times\cA_g$, we have $j^*\overline\alpha_d=0$ for all $d\geq0$ by \cite[Proposition 3.4]{GreerLian2025Level}.
Because $H^1(\overline{\mathcal M}_{1,1},\Q)=0$ and $H^{>0}(\cA_1,\Q)=0$, the K\"unneth decomposition shows that $$\ker j^*=H^2(\overline{\mathcal M}_{1,1},\Q)\otimes H^{2g-2}(\cA_g,\Q),$$ and pushforward along the second projection $p$ restricts to an isomorphism $p_*:\ker j^*\xrightarrow{\sim}H^{2g-2}(\cA_g,\Q)$.
Directly from the definitions, $p_*\overline\alpha_d=\alpha_d$ for $d>0$, while for $d=0$ we have $p_*\overline\alpha_0=(-1)^g\lambda_{g-1}/24=\alpha_0$ since $\deg\mathbb E_1=1/24$ on $\overline{\mathcal M}_{1,1}$.
Via $p_*^{-1}$, the modularity of $\sum_{d\geq0}\alpha_d q^d$ in Theorem \ref{thm:levelone} transfers directly to $\sum_{d\geq0}\overline\alpha_d q^d$.
\end{remark}

\printbibliography
\end{document}